\documentclass[12pt]{amsart}

\theoremstyle{definition}

\theoremstyle{remark}

\numberwithin{equation}{section}

\newcommand{\ra}{\rightarrow}

\newcommand{\ot}{\otimes}
\newcommand{\mtc}{\mathcal}

\newcommand{\Lam}{\Lambda}

\newcommand{\al}{\alpha}
\newcommand{\eps}{\epsilon}
\newcommand{\bn}{\begin}

\numberwithin{equation}{section}

\newtheorem{rem}[equation]{Remark}

\newcommand{\dw}{\downarrow}
\newcommand{\uw}{\uparrow}

\newcommand{\ch}{\chi}
\newcommand{\mtr}{\mathrm}

\begin{document}
\author{Sebastian  Burciu}
\address{Inst.\ of Math.\ ``Simion Stoilow" of the Romanian Academy
P.O. Box 1-764, RO-014700, Bucharest, Romania}

\email{smburciu@syr.edu}

\thanks{The research was supported by CEx05-D11-11/04.10.05.}


\begin{abstract}Let $H$ be a finite dimensional semisimple Hopf algebra over an algebraically closed field of characteristic zero. In this note we give a short proof of the fact that a Hopf subalgebra of $H$ is a depth two subalgebra if and only if it is normal Hopf subalgebra.
\end{abstract}
\title{Depth Two Hopf Subalgebras of Semisimple Hopf algebras}
\maketitle
\section{Introduction}
The finite depth theory has its roots in the classification of $II_1$ subfactors. Particularly, depth two theory has been recently extensively investigated by algebraists. This theory is a type of Galois theory for noncommutative rings and has been intensively studied in \cite{kads},  \cite{KN}, \cite{KS}, \cite{kadtr} and other papers. In the theory, the classical Galois group is replaced by a Hopf algebroid. Results from Hopf Galois extensions had also an important influence in the development of this theory.

Hopf algebroids with a separable base algebra are weak Hopf algebras, and they are used in conformal field theory as well as in other subjects.
In a recent paper \cite{kads}, the author used the theory of weak Hopf algebras to prove prove that a depth two Hopf subalgebra of a semisimple Hopf algebra is a normal Hopf subalgebra. This was done by constructing a map from the semisimple Hopf algebra to a weak Hopf algebra. In this note we give a different proof of this result. The character theory for normal Hopf subalgebras developed in \cite{ker}
and \cite{coset} is used.

All algebras and coalgebras in this paper are defined over an algebraically closed ground field $k$ of characteristic zero. For a vector space $V$ by $|V|$ is denoted the dimension $\mtr{dim}_{k}V$. The comultiplication, counit and antipode of a Hopf algebra are denoted by $\Delta$, $\epsilon$ and
$S$, respectively. For a finite dimensional semisimple Hopf algebra $H$ denote by $\mtr{Irr}(H)$ the set of irreducible characters of $H$. 
All the other notations are those used in \cite{Montg}.
\section{Depth two Hopf subalgebras and normal Hopf subalgebras}
A unital subring B is a right depth two subring of A if there is a split
epimorphism of $A-B$ bimodules from some $A^{n}$ onto $A \ot_B A$ \cite{KS}. A similar condition is imposed for left depth two subring. Finitely generated Hopf-Galois extensions satisfy a stronger form of this condition, where the split epimorphism is replaced by an isomorphism of $A-B$ bimodules.

For a separable algebra $A$ with a separable subalgebra $B$ over a filed the following condition is equivalent to the left depth two condition ( see Theorem 2.1, item 6 of \cite{kadtr}):

As a natural transformation between functors from the category of $B$-modules into the category of right $A$-modules, there is a natural monic from $\mtr{Ind}_B^A\mtr{Res}^B_A\mtr{Ind}_B^A$ into $N\mtr{Ind}_B^A$ for some positive integer $N$. In particular, for each pair of simple modules $V_B$ and $W_A$, the number of isomorphic copies of $W$,
\bn{equation}\label{cond}
<\mtr{Ind}_B^A\mtr{Res}^B_A\mtr{Ind}_B^AV,\;W>\leq N<\mtr{Ind}_B^AV,\;W>
\end{equation}
Now suppose $A=H$ is a finite dimensional semisimple Hopf algebra and $B=K$ is a Hopf subalgebra. Since both algebras are semisimple they are also separable.
 

In terms of character theory the above condition is equivalent to the existence of a positive integer $N$ such that
\bn{equation}\label{condch}
m_{ _H}(\al\uw_{ _K}^H\dw^H_{ _K}\uw_{ _K}^H,\;\ch)\leq Nm_{ _H}(\al\uw_{ _K}^H,\;\ch)
\end{equation}
for all irreducible characters $\al \in \mtr{Irr}(K)$ and $\ch \in \mtr{Irr}(H)$. Since both algebras are semisimple left and right depth two extensions coincide in this situation. 
\subsection{Equivalence for depth two condition}\label{eq}
It is easy to see that $\ch$ is a constituent of $\ch\dw^K_{ _H}\uw^K_{ _H}$ for any irreducible $H$-character $\ch \in \mtr{Irr}(H)$.
Then the depth two condition \ref{condch} is equivalent to the fact that $\al\uw_{ _K}^H\dw^K_{ _H}\uw^K_{ _H}$ and $\al\uw^K_{ _H}$ have the same simple $H$-constituents for any irreducible character $\al \in \mtr{Irr}(K)$.
\subsection{Normal Hopf subalgebras}
Let $H$ be a finite dimensional\\ semisimple Hopf algebra over $k$. Then $H$ is also cosemisimple \cite{Lard}. We use the notation $\Lam_{ _H} \in H$ for the idempotent integral of $H$, $ (\eps(\Lam_{ _H})=1)$, and $t_H \in H^*$ for the integral of $H^*$  with $t_H(1)=|H|$. From Proposition 4.1 of \cite{Lar} it follows that the regular character of $H$ is $t_{ _H}$ and therefore
\begin{equation}
\label{f1}t_{ _H}=\sum_{\ch \in \mtr{Irr}(H)}\chi(1)\chi. 
\end{equation}
If $K$ is a Hopf subalgebra of $H$ then $K$ is a semisimple and cosemisimple Hopf algebra \cite{Montg}. A Hopf subalgebra $K$ of $H$ is called normal if $h_1xS(h_2)\in K$ and $ S(h_1)xh_2 \in K$ for all $x \in K$ and $h \in H$. If $H$ is semisimple Hopf algebra as above then $S^2=\mtr{Id}$ (see \cite{Lard}) and $K$ is normal in $H$ if and only if $h_1xS(h_2)\in K$ for all $x \in K$ and $h \in H$.  If $K^+=Ker(\eps)\cap K$ then $K$ is normal Hopf subalgebra of $H$ if and only if $HK^+=K^+H$.  In this situation $H//K:=H/HK^+$ is a quotient Hopf algebra of $H$ via the canonical map $\pi:H\ra H//K$ (see Lemma 3.4.2 of \cite{Montg}). In our settings $K$ is normal in $H$ if and only if $\Lam_{ _K}$ is central in $H$ . For one implication see Lemma 2.16 of \cite{Mas'}. For the other implication, if $\Lam_{ _K}$ is central in $H$ then $HK^+=H(1- \Lam_{_ K})=(1- \Lam_{_ K})H=K^+H$. For a different proof of this fact see Proposition 2.3. of \cite{kads}.
\subsection{Restriction to normal Hopf subalgebras}
Let $H$ be a semisimple Hopf algebra over an algebraically closed field $k$ and let  $K$ be a normal Hopf subalgebra of $H$. Define an equivalence relation on the set  $\mtr{Irr}(H)$  by $\ch \sim \mu$ if and only $m_{ _K}(\ch \dw_K^H,\;\mu\dw_K^H)>0$. This is the equivalence relation $r^{H^*}_{ _{L^*,\;k}}$ from \cite{coset}. It is proven that $\ch \sim \mu$ if and only if $\frac{\ch\dw_K^H}{\ch(1)}=\frac{\mu\dw_K^H}{\mu(1)}$ (see Theorem 4.3 of \cite{coset}). Thus the restriction of $\ch$ and $\mu$ to $K$ either have the same irreducible constituents or they don't have common constituents at all.

The above equivalence relation determines an equivalence relation on the set of irreducible characters of $K$. Two irreducible $K$-characters $\al $ and $\beta$ are equivalent if and only if they are constituents of  $\ch\dw_{ _K}^H$ for some irreducible character $\ch$ of $H$. Let $\mtc{C}_1,\cdots \mtc{C}_s$ be the equivalence classes of the above equivalence relation on $\mtr{Irr}(H)$. Let $\mtc{A}_1,\cdots \mtc{A}_s$ be the corresponding equivalence classes on $\mtr{Irr}(K)$. The formulae from Section of $4.1$ of \cite{coset} imply that if $\ch \in \mtc{C}_i$ then
\bn{equation} \label{restrform}
\ch\dw_{ _K}^H=\frac{\ch(1)}{|\mtc{A}_i|}\sum_{\al\in \mtc{A}_i}\al(1)\al
\end{equation}
 where $|\mtc{A}_i|=\sum_{\al\in \mtc{A}_i}\al(1)^2$.
Also if $\al \in \mtc{A}_i$ then
\bn{equation}\label{indform}
\al\uw_{ _K}^H=\frac{\al(1)}{a_i(1)}\frac{|H|}{|K|}\sum_{\ch \in \mtc{C}_i}\ch(1)\ch
\end{equation}
where $a_i(1)=\sum_{ \ch \in \mtc{C}_i}\ch(1)^2$.
\begin{rem}\label{rch}
These two formulae show that $\al\uw_{ _K}^H\dw_{ _K}^H\uw_{ _K}^H$ has the same irreducible $K$-constituents as $\al\uw_{ _K}^H$ for all irreducible characters $\al \in \mtr{Irr}(K)$. As we have seen in subsection \ref{eq} this implies that $K$ is a depth two subalgebra of $H$.
\end{rem}

\section{Proof of the main result}
Let $\eps_{ _K}$ be the character of the trivial $K$-module.
The following lemma will be used in the proof of the main theorem. This is a slightly weakened version of Corollary 2.5 from \cite{ker}.
\bn{lemma}
Let $H$ be a finite dimensional semisimple Hopf algebra and $K$ be a Hopf subalgebra of $H$. Then $K$ is normal if and only if $\eps_{ _K}\uw^H_K\dw_K^H=\frac{|H|}{|K|}\eps_{ _K}$.
\end{lemma}
\bn{proof}
Suppose $K$ is normal Hopf subalgebra of $H$. Then Corollary 2.5 of \cite{ker} implies the desired condition.

Conversely, suppose that the above condition is satisfied. By Frobenius reciprocity $m_{ _K}(\ch \dw_K, \; \eps_{ _K})=m_{ _H}(\ch,\; \eps\uparrow_{ _K}^H)$.
The above condition implies that for any irreducible character $\chi$ of $H$ we have that the value of $m_{ _K}(\ch \dw_K, \; \eps_{ _K})$ is either $\ch(1)$ if $\ch$ is a
constituent of $\eps \uw_{ _K}^H$ or $0$ otherwise. But if $\Lam_{ _K}$ is
the idempotent integral of $K$ then $m_{ _K}(\ch \dw_{ _K}, \;
\eps_{ _K})=\ch(\Lambda_{ _K})$ (see Proposition 1.7.2 of \cite{N}). Thus $\ch(\Lam_{ _K})$ is either zero or $\ch(1)$ for any irreducible character $\ch$ of $H$. This implies that
$\Lambda_{ _K}$ is a central idempotent of $H$ and therefore $K$ is a normal
Hopf subalgebra of $H$.
\end{proof}
\bn{thm} Let $H$ be a finite dimensional semisimple Hopf algebra. A Hopf subalgebra $K$ of $H$ is depth two subalgebra if and only if $K$ is normal in $H$.
\end{thm}
\bn{proof} 
If $K$ is a normal Hopf algebra then Remark \ref{rch} shows that $K$ is a depth two subalgebra.

Suppose now that $K$ is a depth two subalgebra of $H$. Using Frobenius reciprocity the condition \ref{condch} is equivalent to:
\bn{equation}\label{cond2}
m_{ _K}( \al\uw_{ _K}^H\dw_{ _K}^H,\;\ch\dw_{ _K}^H)\leq Nm_{ _K}(\al,\;\ch\dw_{ _K}^H)
\end{equation}
for all $\al \in \mtr{Irr}(K) $ and $\ch \in \mtr{Irr}(H)$.
If $\ch= \eps_{ _H}$ is the trivial $H$-character then its restriction to $K$ is the trivial $K$-character and the above condition shows that $m_{ _K}(\al\uw_K^H\dw_K^H,\;\eps_{ _K})=0$ if $\al\neq \eps_{ _K}. $

The regular character of $K$ induced to $H$ and then restricted back to $K$ is  the regular character of $K$ multiplied $\frac{|H|}{|K|}$. Thus $$t_{ _K}\uw_K^H\dw_K^H=\frac{|H|}{|K|}t_{ _K}.$$
On the other hand using formula \ref{f1} for $K$ one has $$t_{ _K}\uw_K^H\dw_K^H=\sum_{\al \in \mtr{Irr}(K)}\al(1)\al\uw_K^H\dw_K^H.$$ Thus the multiplicity of $\eps_{ _K}$ in the above expression is $\frac{|H|}{|K|}$. But in the above sum, $\eps_{ _K}$ might be constituent only in the term corresponding to the trivial $K$-character $\al=\eps_{ _K}$ since $m_{ _K}(\al\uw_K^H\dw_K^H,\;\eps_{ _K})=0$ for $\al\neq \eps_{ _K}. $
Therefore $m_{ _K}(\eps_{ _K}\uw_K^H\dw_K^H,\;\eps_{ _K})=\frac{|H|}{|K|}\eps_K$ and a dimension argument implies that  $\eps_{ _K}\uw^H_{ _K}\dw_{ _K}^H=\frac{|H|}{|K|}\eps_{ _K}$ . The previous lemma shows that $K$ is normal in $H$.

\end{proof}
\bibliographystyle{plain}
\bibliography{depthtwo}
\end{document}